\documentclass[12pt]{amsart}
\usepackage{amssymb, amsmath, amsthm,amsrefs,marvosym,wasysym}
\usepackage{latexsym}
\usepackage{color}
\usepackage{exscale}
\usepackage{fullpage}
\usepackage{rotating}
\usepackage{bm, bbm, mathrsfs}
\usepackage{graphics, graphicx}
\newtheorem{theorem}{Theorem}[section]

\newtheorem{claim}[theorem]{Claim}

\newtheorem{definition}{Definition}

\newtheorem{assumption}{Assumption}
\newtheorem{remark}{Remark}[section]

\def\vp{\varphi}

\def\eq#1{(\ref{#1})}
\def\nn{\nonumber}
\def\({\left(\begin{array}{cccccc}}
\def\){\end{array}\right)}

\def\eq#1{(\ref{#1})}
\def\nn{\nonumber}

\def\({\left(\begin{array}{cccccc}}
\def\){\end{array}\right)}

\def\bes{\begin{eqnarray}}
\def\ees{\end{eqnarray}}

\newcommand{\beq}{\begin{equation}}
\newcommand{\eeq}{\end{equation}}
\newcommand{\bea}{\begin{eqnarray}}
\newcommand{\eea}{\end{eqnarray}}
\newcommand{\beann}{\begin{eqnarray*}}
\newcommand{\eeann}{\end{eqnarray*}}

\newcommand{\RR}{\mathbb{R}}

\newcommand{\br}{\bar r}

\DeclareMathOperator{\supp}{supp}

\DeclareMathOperator{\sgn}{sgn}

\newcommand{\s}{\ensuremath{\mathrm{s}}}

\newcommand{\w}{\ensuremath{\mathrm{w}}}

\newcommand{\kaa}{\ensuremath{\mathrm{k}}}

\newcommand{\lj}{\big[\!\!\big[}
\newcommand{\rj}{\big]\!\!\big]}

\DeclareMathOperator{\grad}{grad}
\DeclareMathOperator{\dv}{div}

\numberwithin{equation}{section}
 
\begin{document}

\title{Multi-d Isothermal Euler Flow: Existence of unbounded radial similarity solutions}

\begin{abstract}
	We show that the multi-dimensional compressible Euler system for isothermal 
	flow of an ideal, polytropic gas admits global-in-time, radially symmetric solutions with unbounded 
	amplitudes due to wave focusing. The examples are similarity solutions and involve 
	a converging wave focusing at the origin.  At time of collapse, the density, but not 
	the velocity, becomes unbounded, resulting in an expanding shock wave. 
	The solutions are constructed as functions of radial distance to the 
	origin $r$ and time $t$. We verify that they provide genuine, weak 
	solutions to the original, multi-d, isothermal Euler system.
	
	While motivated by the well-known Guderley solutions to the full Euler 
	system for an ideal gas, the solutions we consider are of a different type.
	In Guderley solutions an incoming shock propagates toward the origin 
	by penetrating a stationary and ``cold'' gas at zero pressure (there is no counter 
	pressure due to vanishing temperature near the origin),
	accompanied by blowup of velocity and pressure, but not of density, at collapse. 
	It is currently not known whether the full system admits unbounded solutions in 
	the absence of zero-pressure regions. The present work shows that the 
	simplified isothermal model does admit such behavior. 
\end{abstract}

\author{Helge Kristian Jenssen }\address{H.~K.~Jenssen, Department of
Mathematics, Penn State University,
University Park, State College, PA 16802, USA ({\tt
jenssen@math.psu.edu}).}

\author{Charis Tsikkou}\address{C. Tsikkou, Department of
Mathematics, West Virginia University,
Morgantown, WV 26506, USA ({\tt
tsikkou@math.wvu.edu}).}

\date{\today}
\maketitle

\tableofcontents

%%%%%%%%%%%%%%%%%%%%%%%%%%%%%%%%%%%%%%%%%%%%%%%%
%%%%%%%%%%%%%%%%%%%%%%%%%%%%%%%%%%%%%%%%%%%%%%%%
\section{Equations}
%%%%%%%%%%%%%%%%%%%%%%%%%%%%%%%%%%%%%%%%%%%%
The compressible Euler system for barotropic flow in $\RR^n_{\bf  x}$ is given by
\begin{align}
	\rho_t+\dv_{\bf x}(\rho \bf u)&=0 \label{mass_m_d_isthrml_eul}\\
	(\rho{\bf  u})_t+\dv_x[\rho {\bf  u}\otimes{\bf  u}]+\grad_{\bf  x} p&=0,\label{mom_m_d_isthrml_eul}
\end{align}
where the independent variables are time $t$ and position ${\bf  x}$, and the primary dependent variables are 
density $\rho$ and velocity ${\bf  u}$, while pressure is a given function of density, $p=p(\rho)$.
In {\em isothermal} flow of an ideal, polytropic gas the pressure is a linear function of density:
\beq\label{pressure}
	p(\rho)=a^2\rho\qquad\qquad\text{($a>0$ constant)}.
\eeq
For {\em radial} ($\equiv$ spherically symmetric) solutions the dependent variables are functions 
of time $t$ and radial distance $r=|{\bf x}|$ to the origin, and the velocity field is purely radial: 
${\bf u}=u\frac{\bf x}{r}$. In this case \eq{mass_m_d_isthrml_eul}-\eq{mom_m_d_isthrml_eul} 
reduces to a quasi one-dimensional system:
\begin{align}
	\left(r^m\rho\right)_t+\left(r^m\rho u\right)_r &= 0\label{mass}\\
	\left(r^m\rho u \right)_t+\left(r^m(\rho u^2+p)\right)_r &= mr^{m-1}p,\label{momentum}
\end{align}
where $m=n-1$. For smooth (Lipschitz) flows this reduces further to
\begin{align}
	\rho_t+u\rho_r+\rho\left(u_r+\frac{mu}{r}\right) &= 0\label{m_eul}\\
	u_t+ uu_r +\frac{p_r}{\rho}&= 0.\label{mom_eul}
\end{align}
In this work we shall be concerned exclusively with complete radial isothermal 
flows of {\em similarity type}. This means $\rho(t,r)$ and $u(t,r)$ are defined for 
all $t\in\RR$, $r>0$, and are of the form 
\beq\label{fncl_relns}
	\rho(t,r)=\sgn(t)|t|^\beta\Omega(\xi)\,,\qquad u(t,r)=U(\xi),
\eeq
where the {\em similarity variable} $\xi$ is given by
\[\xi=\frac{r}{t}.\]
A discussion of our results and their relations to earlier works appears
in Section \ref{conv_div_isothermal}. 

At this stage $\beta\in\RR$ in \eq{fncl_relns} is a free parameter. Substitution of \eq{pressure} and 
\eq{fncl_relns} into \eq{m_eul}-\eq{mom_eul} yields the similarity ODEs 
(where ${}'\equiv \frac{d}{d\xi}$)
\begin{align}
	(U-\xi)\frac{\Omega'}{\Omega}+U'+\Big(\beta+\frac{mU}{\xi}\Big)&=0
	\label{isoth_omega_ode}\\
	a^2\frac{\Omega'}{\Omega}+(U-\xi)U'&=0.
	\label{isoth_u_ode}
\end{align}
Solving for $\frac{\Omega'}{\Omega}$ in \eq{isoth_u_ode} and substituting 
into \eq{isoth_omega_ode} yield a single ODE for $U$:
\beq\label{U_ode}
	U'=\frac{a^2}{(U-\xi)^2-a^2}\Big(\beta+\frac{mU}{\xi}\Big).
\eeq
Using this in \eq{isoth_u_ode} gives
\beq\label{Omega_ode}
	\frac{\Omega'}{\Omega}=-\frac{U-\xi}{(U-\xi)^2-a^2}\Big(\beta+\frac{mU}{\xi}\Big).
\eeq
Before analyzing the similarity ODEs we consider the jump conditions
in similarity variables.

%%%%%%%%%%%%%%%%%%%%%%%%%%%%%%%%%%%%%%%%%%%%
%%%%%%%%%%%%%%%%%%%%%%%%%%%%%%%%%%%%%%%%%%%%
\section{Rankine-Hugoniot and Entropy conditions for similarity flows}
%%%%%%%%%%%%%%%%%%%%%%%%%%%%%%%%%%%%%%%%%%%%
Consider the radial barotropic Euler system \eq{mass}-\eq{momentum}, and 
assume that a discontinuity propagates along the path $r=\mathcal R(t)$. 
The Rankine-Hugoniot conditions are then
\begin{align}
	\dot{\mathcal R}\lj\rho\rj &= \lj \rho u\rj \label{rh_1}\\
	\dot{\mathcal R}\lj\rho u\rj &= \lj \rho u^2+p\rj, \label{rh_2}
\end{align}
where $\dot{} \equiv \frac{d}{dt} $.
Here and below we use the convention that, for any quantity $q=q(t,r)$,
$\lj q\rj$ denotes the jump in $q$ as $r$ decreases, i.e.,
\[\lj q\rj:=q_+-q_-\equiv q(t,\mathcal R(t)+)-q(t,\mathcal R(t)-).\]
Next, denoting the local sound speed by
\[c:=\sqrt{p'(\rho)},\]
the entropy condition for a 1-shock requires that 
\beq\label{e1}
	u_--c_-> \dot{\mathcal R}> u_+-c_+,
\eeq
while the entropy condition for a 2-shock requires that 
\beq\label{e2}
	u_-+c_-> \dot{\mathcal R}> u_++c_+.
\eeq

%%%%%%%%%%%%%%%%%%%%%%%%%%%%%%%%%%%%%%%%%%%%
\subsection{Radial isothermal similarity shocks}
We next specialize to ``similarity shocks'' in radial isothermal flow: the pressure law is given 
by \eq{pressure} and the shock is assumed to propagate along a path of the form 
$\xi\equiv \bar \xi$, i.e., $\mathcal R(t)=\bar\xi t$. Furthermore, it is assumed  
that the density and velocity on either side of the shock are of the form \eq{fncl_relns},
with $\beta$ taking the same value on both sides. 
Let  $(U_+,\Omega_+)$ and $(U_-,\Omega_-)$ denote the parts of the solution on the 
outside and inside of the shock, respectively. (``Outside'' and ``inside'' refer to further 
away from and closer to $r=0$, respectively.)

The Rankine-Hugoiniot conditions reduce to 
\begin{align*}
	\bar\xi\lj\Omega\rj &= \lj \Omega U\rj \\
	\bar\xi\lj\Omega U\rj &= \lj \Omega (U^2+a^2)\rj, 
\end{align*}
where $\lj\cdot\rj$ now denotes jump across $\xi=\bar\xi$. The entropy conditions 
\eq{e1}-\eq{e2} take the form
\begin{align}
	&U_-(\bar\xi)> \bar\xi+a> U_+(\bar\xi)\qquad\text{for a 1-shock}\label{isoth_e_1}\\
	&U_-(\bar\xi)> \bar\xi-a> U_+(\bar\xi)\qquad\text{for a 2-shock.}\label{isoth_e_2}
\end{align}
In particular, these relations show that for any shock in radial isothermal flow, the velocity 
necessarily decreases as we traverse the shock from the inside to the outside.
%We also note that, for isothermal flow, the entropy conditions \eq{e1}-\eq{e2} involve 
%only the velocity $u$ (and not the density field $\rho$; this feature is particular to 
%isothermal flow.)

Finally, setting $V_\pm:=U_\pm-\bar\xi$, where $U_\pm$ denotes $U_\pm(\bar\xi)$,
the Rankine-Hugoniot conditions take the form $\lj\Omega V\rj=0$ and $\lj\Omega VU+a^2\Omega\rj=0$.
It follows from these that $V_+V_-=a^2$, and that
\beq\label{+ito-}
	U_+=\bar\xi+\frac{a^2}{U_--\bar\xi}\qquad \text{and}\qquad\Omega_+=\frac{(U_--\bar\xi)^2}{a^2}\Omega_-.
\eeq
Alternatively, solving for $V_-$ and $\Omega_-$, we have
\beq\label{-ito+}
	U_-=\bar\xi+\frac{a^2}{U_+-\bar\xi}\qquad \text{and}\qquad\Omega_-=\frac{(U_+-\bar\xi)^2}{a^2}\Omega_+.
\eeq

%%%%%%%%%%%%%%%%%%%%%%%%%%%%%%%%%%%%%%%%%%%%
%%%%%%%%%%%%%%%%%%%%%%%%%%%%%%%%%%%%%%%%%%%%
\section{Converging-diverging isothermal flows}\label{conv_div_isothermal}
%%%%%%%%%%%%%%%%%%%%%%%%%%%%%%%%%%%%%%%%%%%%
By a ``converging-diverging solution'' we shall mean a radial similarity solution in which 
a wave approaches the origin, ``collapses'' there at some instant in time, resulting  
in a reflected wave moving away from the origin. Without loss of generality we set 
the time of collapse to be $t=0$. 

We shall search for this type of solutions within the class of isothermal similarity solutions 
introduced above. To be of physical interest the solutions should satisfy, 
as a minimum, the following requirements:
\begin{itemize}
	\item[(A)] the velocity vanishes along $\{r=0\}$: $u(t,0)\equiv 0$;  
	\item[(B)] at any fixed location $r>0$, the limits
	\[\lim_{t\to 0}u(t,r)\qquad\text{and}\qquad \lim_{t\to 0}\rho(t,r)\]
	both exist as finite numbers. (Note that this requirement leaves open the possibility 
	that $\rho(0,r)$ and/or $u(0,r)$ may blow as $r\downarrow 0$.)
\end{itemize}
In addition we shall require that the density field is everywhere strictly positive:
\begin{itemize}
	\item[(C)] the density never vanishes: $\rho(t,r)> 0$ for all $t\in\RR$, $r\geq 0$.
\end{itemize}
Further constraints will be imposed later to guarantee that the solutions, as function of 
$(t,{\bf x})\in\RR\times \RR^n$, provide genuine weak solutions of the original, multi-d isothermal 
system \eq{mass_m_d_isthrml_eul}-\eq{mom_m_d_isthrml_eul}. In particular, we shall 
require that the conserved quantities map time continuously into $L^1_{loc}(\RR^n)$;
see Section \ref{weak_solns} and also Section \ref{final_rmks}.

For the {\em full} Euler system (including conservation of energy) the seminal work 
\cites{gud} by Guderley established the existence of converging-diverging similarity solutions
in which a shock wave propagates into a quiescent state near the
origin, focuses (collapses) at the origin, and reflects an expanding shock wave. 
Building on the detailed work of Lazarus \cite{laz} (which also treats the case of a 
collapsing vacuum), the present authors recently showed in \cite{jt1} that these ``Guderley solutions'' 
provide examples of genuine, entropy admissible, weak solutions to the full, multi-d Euler system. 
A key feature of these converging-diverging shock solutions 
is that they provide concrete Euler flows suffering pointwise blowup 
of primary flow variables (as opposed to blowup of their gradients).

Although the Guderley solutions establish the possibility of amplitude blowup in Euler 
flows for ideal gases, they are also at the borderline
of the regime where one would expect the Euler system to be physically accurate. More
precisely, in order to provide an {\em exact} weak solution, the sound speed in the 
quiescent state that the incoming shock moves into must vanish. For the 
ideal gas case under consideration, this means that the incoming shock does not 
experience any upstream counter-pressure. (The gas is at zero temperature there
and this is sometimes referred to as a  ``cold gas assumption.'') 
It appears reasonable that this lack of 
counter-pressure facilitates unbounded growth of the shock speed, with concomitant
increases in pressure and temperature. It is unclear at present whether this is the 
(or part of the) mechanism driving the blowup in Guderley solutions for the full Euler system. 
The alternative is that the blowup is a purely geometric effect driven by {\em wave 
focusing}, much like what occurs for radial solutions of the linear, multi-d wave equation.

{\em The main goal of the present work is to show that amplitude blowup can 
occur in converging-diverging flows for the simplified isothermal Euler model, 
even in the presence of an everywhere strictly positive pressure field.} To the best of our 
knowledge, the solutions we generate are the first examples of unbounded barotropic flows
that meet the requirements (A)-(C) above.
While these isothermal solutions are qualitatively different from the Guderley solutions 
for the full system described earlier (in particular, they are continuous up to collapse), 
they indicate that the real agent for blowup is the focusing of waves 
at the center of motion. On the other hand, it still remains an open problem to 
exhibit concrete flows for the full Euler system that exhibit blowup in the absence of
zero-pressure regions.

For completeness we include some remarks on what 
is known about radial Euler flows with ``general'' initial data. 
First, there is currently no result for the full, multi-d 
Euler system, radial or not, that guarantees global-in-time existence.
For radial isentropic flows, i.e., solutions to \eq{mass}-\eq{momentum} 
with $p(\rho)=a^2\rho^\gamma$ and $\gamma>1$, 
results by Chen-Perepelitsa \cite{cp} and Chen-Schrecker \cite{cs} 
provide existence of weak, finite energy 
solutions via the method of compensated compactness. In fact,
the recent work \cite{sch} is the first to show that the solutions one
obtains in this manner provide genuine, weak solutions to the original, 
multi-d isentropic Euler system 
\eq{mass_m_d_isthrml_eul}-\eq{mom_m_d_isthrml_eul} 
on {\em all} of space. On the other hand, there appears to be little 
hope of extending this approach (i.e., compensated compactness)
to the radial full system, or even (for technical reasons \cite{sch1}) 
to the radial, isothermal ($\gamma=1$) system. 

As far as we know, the currently strongest, global existence result for the radial
isothermal system applies to the case of {\em external} flows, i.e., 
for flows outside of a fixed ball. This problem was analyzed in \cite{mmu1}
by exploiting the Glimm scheme, providing existence for a certain class
of initial data of bounded variation; for an extension, see \cite{mmu2}. 
The results of the present paper shows that, in order to extend these results to 
solutions defined on {\em all} of space (i.e., including the origin), one must necessarily 
contend with unbounded solutions. 
 
For results closer to the present work, which concerns concrete 
Euler flows in several space dimensions, see Chapter 7 of Zheng's 
monograph \cite{zheng} on multi-d Riemann problems, some of which 
generate purely radial flows. However, we stress that the radial flows we 
construct below are not solutions to Riemann problems. Specifically, 
the solutions we display are necessarily non-constant in the radial 
direction at all times.

The rest of the present paper is organized as follows. Section \ref{constr_conv_div_solns}
provides a detailed construction of the radial speed $u(t,r)$ and the corresponding density 
$\rho(t,r)$ for converging-diverging similarity flows for the isothermal Euler system.
In Section \ref{weak_solns} we briefly recall the definition of weak solutions to the barotropic 
Euler system, including its formulation for the special case of radial solutions.  
In Section \ref{sim_weak_solns} we verify that the radial similarity flows we
construct provide genuine weak solutions to the original, multi-d isothermal Euler system.
The main result is summarized in Theorem \ref{main_result}. Finally, Section \ref{final_rmks}
collects some additional observations about the flows constructed in this paper.

%%%%%%%%%%%%%%%%%%%%%%%%%%%%%%%%%%%%%%%%%%%%
%%%%%%%%%%%%%%%%%%%%%%%%%%%%%%%%%%%%%%%%%%%%
\section{Construction of converging-diverging isothermal flows}\label{constr_conv_div_solns}
%%%%%%%%%%%%%%%%%%%%%%%%%%%%%%%%%%%%%%%%%%%%
To construct concrete examples of converging-diverging isothermal similarity flows, we 
start with the ODE \eq{U_ode} for the velocity $U(\xi)$. This ODE has three 
critical points: the origin $(0,0)$ and the points $\pm P_{\w}:=(\pm\xi_{\w},\pm U_{\w})$, where
\[\xi_{\w}:=-\frac{am}{m+\beta},\qquad U_{\w}:=\frac{a\beta}{m+\beta}.\] 
(The subscript ``w'' stands for ``weak,'' for reasons to be clear later.) 
We also observe that its solutions are symmetric about the origin: if $\xi\mapsto U(\xi)$ 
is a solution of \eq{U_ode}, so is $\xi\mapsto -U(-\xi)$.

Instead of performing a lengthy analysis of all possible cases, from now on we focus 
on the cases where
\beq\label{case3}
	-m<\beta<0 \qquad\text{and}\qquad m=1\quad\text{or}\quad m=2.
\eeq
In particular, $\xi_{\w}<0$ and  $U_{\w}<0$ for all cases under consideration. Introducing the straight lines
\[l_\pm:=\{U=\xi\pm a\}\qquad\text{and}\qquad \omega:=\{\beta+\textstyle\frac{mU}{\xi}=0\},\]
we have that $\pm P_{\w}=l_\pm\cap\omega$. Linearizing \eq{U_ode} about the 
critical points $\pm P_\w$, we set
\beq\label{lambdas}
	\lambda_\pm=\textstyle\frac{1}{2}\Big[\big(1+\frac{m}{2}(1+\mu)\big)
	\pm\sqrt{\big(1+\frac{m}{2}(1+\mu)\big)^2-2m(1+\mu)^2}\Big],
\eeq
where 
\[\mu:=\textstyle\frac{\beta}{m}\in(-1,0).\]
It is immediate to verify that the radicand in \eq{lambdas} is strictly positive 
whenever \eq{case3} holds. An analysis of the critical points shows that:
\begin{enumerate}
	\item[(a)] The point $P_{\w}$ is an unstable node for \eq{U_ode} whenever \eq{case3} holds, i.e., 
	we have $0<\lambda_-<\lambda_+$.
	\item[(b)] There are two solutions leaving $P_{\w}$ along the directions $\pm(1,1-\lambda_+)$.
	\item[(c)] All other solutions leaving $P_{\w}$ do so along the directions $\pm(1,1-\lambda_-)$.
	\item[(d)] Whenever \eq{case3} holds we have $1-\lambda_+<0$ and $-\mu<1-\lambda_-<1$;
	thus all but the two solutions described in (b), enter the region between the straight lines
	$\omega$ and $l_+$.
	\item[(e)] There is a unique solution passing through $(0,0)$; it does so with slope $-\frac{\beta}{n}$,
	and this solution is located below $l_+$ and above $\omega$; it extends back 
	(i.e., as $\xi$ decreases) to $P_{\w}$, approaching $P_{\w}$ along the direction $-(1,1-\lambda_-)$.
\end{enumerate}
We denote the unique solution described in (e) by $\hat U(\xi)$. It passes through the origin 
and, by symmetry about the origin, is defined for all $\xi\in[\xi_{\w},-\xi_{\w}]$,
and connects to the third critical point $-P_\w$.
See Figure 1.

%%%%%%%%%%%%%%%%%%
%	FIGURE 
%%%%%%%%%%%%%%%%%%
\begin{figure}\label{Figure_1}
	\centering
	\includegraphics[width=16cm,height=8cm]{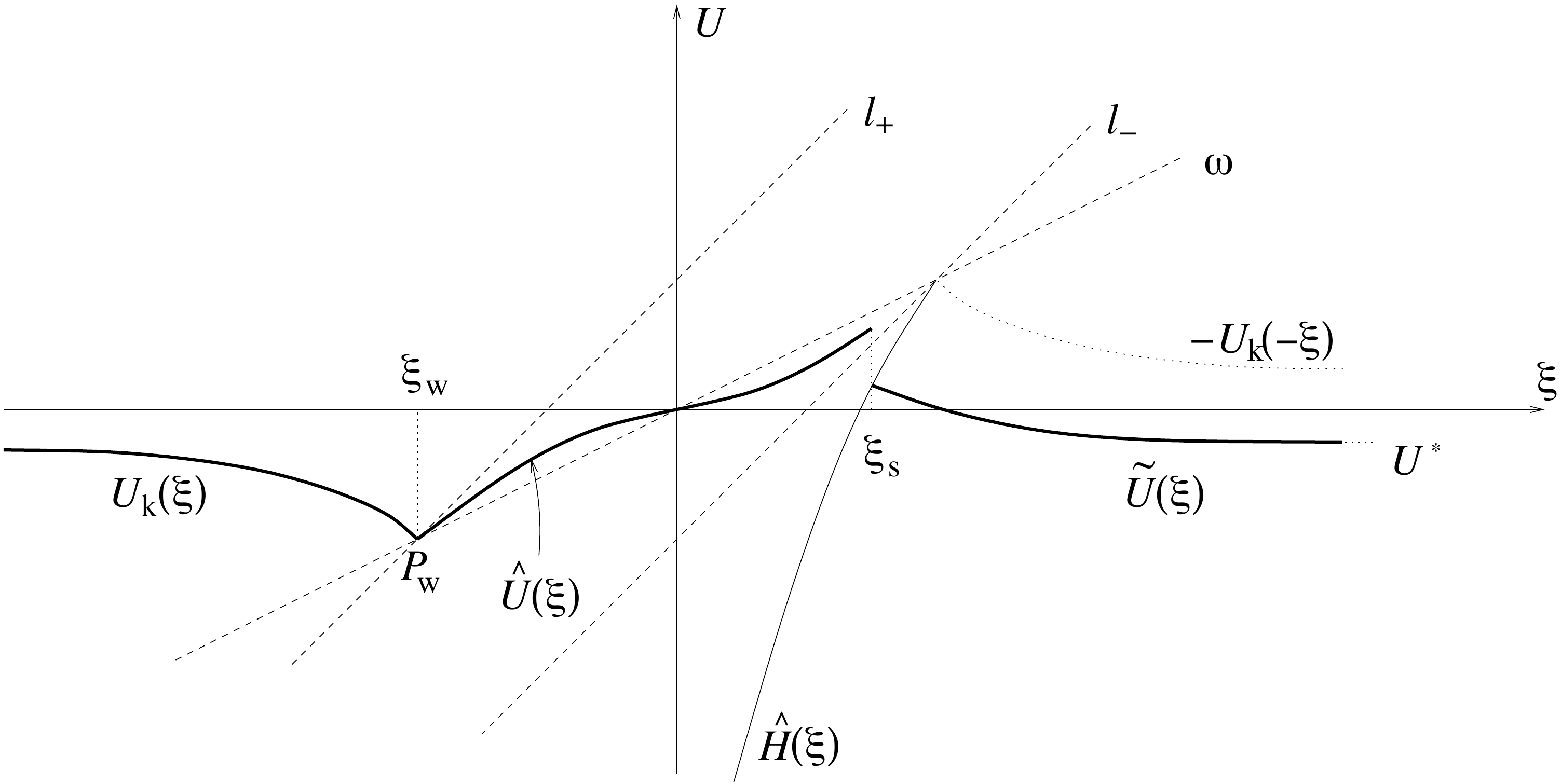}
	\caption{Complete $U(\xi)$-profile (schematic).}
\end{figure}

%%%%%%%%%%%%%%%%%%%%%%%%%%%%%%%%%%%%%%%%%%%%
\subsection{The radial speed $u(t,r)$ for $t\leq 0$}
%%%%%%%%%%%%%%%%%%%%%%%%%%%%%%%%%%%%%%%%%%%%
The part of $\hat U(\xi)$ corresponding to $\xi\in[\xi_{\w},0]$ yields, via \eq{fncl_relns}${}_2$,
the radial speed $u(t,r)$ within the sector 
\[S_-:=\{(r,t)\,:\,\xi_{\w}\leq\textstyle \frac{r}{t}\leq 0\}\] 
in the $(r,t)$-plane. Note that the choice of the solution $\hat U(\xi)$ in this region is dictated by 
requirement (A) above. Similarly, we shall use a certain portion of $\hat U(\xi)$ for $\xi>0$ 
to obtain the radial speed $u(t,r)$ within a sector 
\[S_+:=\{(r,t)\,:\,0\leq \textstyle \frac{r}{t}\leq \xi_\s\}.\] 
Here the value of $\xi_\s\in(0,-\xi_\w)$, yet to be determined, corresponds to the path $t\mapsto \xi_\s t$
of an expanding shock wave for $t>0$.

However, we first need to continue the relevant $U$-solution beyond $\xi_{\w}$, all the way  
down to $\xi=-\infty$.
Now, there are infinitely many solutions of \eq{U_ode} defined for all $\xi<\xi_{\w}$, passing 
through $P_{\w}$, and with the property that they enter (as $\xi$ decreases) the region $\mathcal U$ to the left of
$P_\w$ and above $\omega$, i.e.,
\[\mathcal U:=\{\,(\xi,U)\,:\, \xi< \xi_\w\quad \text{and}\quad U>-\mu\xi\,\}.\]
Let $\check U(\xi)$ denote any such solution.
We therefore have an infinity of choices for $\check U(\xi)$. 
As we shall see below, all of these solutions (that enter $\mathcal U$ at points along $\omega$)
tend to finite limits at $\xi=-\infty$, as dictated by the first part of requirement (B) above. 
However, it will be convenient for the subsequent analysis to also have $\check U(-\infty)<0$.
We proceed to show that there are solutions satisfying this constraint, as well as the 
constraints in \eq{case3}.

%%%%%%%%%%%%%%%%%%%%%%%%%%%%%%%%%%%%%%%%%%
\subsubsection{Asymptotics for large, negative $\xi$-values.}
%%%%%%%%%%%%%%%%%%%%%%%%%%%%%%%%%%%%%%%%%%
As is clear from the linearization of \eq{U_ode} at $P_\w$, all but one of the 
solutions $\check U(\xi)$ defined on $(-\infty,\xi_\w)$ approach $P_\w$ along $(1,1-\lambda_-)$;
all of these connect smoothly at $\xi=\xi_\w$ with the solution $\hat U(\xi)$ on $[\xi_w,0]$ 
considered above. 
%(Note that, according to \eq{U_ode}, any solution $\check U(\xi)$ satisfies 
%$\check U'(\xi)\leq0$ within $\mathcal U$; in particular, any solution 
%entering $\mathcal U$ at a point $(\xi_0,-\mu\xi_0)$ on $\omega$, remains within $\mathcal U$
%for all $\xi\leq\xi_0$.)
The exception is the ``kink-solution'' $U_\kaa(\xi)$ which approaches $P_\w$ 
along $(1,1-\lambda_+)$. It is clear that $U_\kaa(\xi)$ lies above any  
solution $\check U(\xi)$ of \eq{U_ode} which is located in $\mathcal U$ and which exits 
$\mathcal U$ (as $\xi$ increases) at a point on $\omega$. 
For our purpose of having $\check U(-\infty)<0$, it therefore suffices to identify cases 
for which $U_\kaa(\xi)$ tends to a strictly negative limit at $\xi=-\infty$, and then employ $U_\kaa$
in our construction of $u(t,r)$ within the sector 
\[S'_-:=\{(r,t)\,:\,-\infty<\textstyle \frac{r}{t}\leq \xi_{\w}\}.\]  
We start by observing that for $(\xi,U)\in\mathcal U$ we have
$U-\xi-a\geq U+\mu \xi\geq 0$ so that
\[\frac{U+\mu \xi}{U-\xi-a}\leq 1\qquad\text{within $\mathcal U$}.\]
Therefore, any solution $\check U(\xi)$ of \eq{U_ode} in $\mathcal U$ satisfies
\[\check U'(\xi)=\frac{a^2m(\check U(\xi)+\mu \xi)}{\xi(\check U(\xi)-\xi-a)(\check U(\xi)-\xi+a)}
\geq \frac{a^2m}{\xi(\check U(\xi)-\xi+a)}.\]
Specializing to the kink-solution $U_\kaa(\xi)$, which satisfies 
$U_\kaa(\xi)>U_\w$ for $\xi<\xi_\w$, we obtain
\[U_\kaa'(\xi) \geq \frac{a^2m}{\xi(U_\kaa(\xi)-\xi+a)}
>\frac{ma^2}{\xi(U_\w+a-\xi)}\qquad\text{for $\xi<\xi_\w$.}\]
Integrating from $\xi=-\infty$ to $\xi=\xi_\w$, and using that $U_\kaa(\xi_\w)=U_\w$, 
yields
\beq\label{limit}
	U_\kaa(-\infty)< U_\w+a^2m\int_{-\infty}^{\xi_\w}\frac{d\xi}{\xi(\xi-(a+U_\w))}.
%	&=a\Big[\frac{\mu}{1+\mu}+m\frac{(1+\mu)}{(1+2\mu)}\log(2(1+\mu))\Big].
\eeq
Therefore, whenever $m$ and $\beta$ satisfy $-m<\beta<0$, and are
such that the right-hand side of \eq{limit} is non-positive, then the 
kink-solution $U_\kaa(\xi)$ tends to a strictly negative limit as $\xi\downarrow -\infty$. 
E.g., with $m=2$ and $\beta=-1$, the right-hand side of \eq{limit} takes the value zero,
while for $m=1$ and $\beta=-\frac{1}{2}$ it takes a strictly negative value.

\begin{assumption}\label{assmpn}
	From now on it is assumed that $m$ and $\beta$ are such that $m=1$ or $m=2$, 
	$-m<\beta<0$, and at the same time 
	\[U^*:=U_\kaa(-\infty)<0;\] 
	the argument above demonstrates that such values of $m$ and $\beta$ exist.
\end{assumption}
\noindent As indicated above, we use the kink-solution $U_\kaa(\xi)$ to specify the 
radial speed $u(t,r)$, via \eq{fncl_relns}${}_2$, within the sector 
$S'_-:=\{(r,t)\,:\,-\infty<\textstyle \frac{r}{t}\leq \xi_{\w}\}.$

%%%%%%%%%%%%%%%%%%%%%%%%%%%%%%%%%%%%%%%%%%%%
\subsection{The radial speed $u(t,r)$ for $t\geq 0$; the reflected shock.}\label{u_for_pos_t}
%%%%%%%%%%%%%%%%%%%%%%%%%%%%%%%%%%%%%%%%%%%%
Next, we need to specify the radial speed $u(t,r)$
within the sector 
\[S'_+:=\{(r,t)\,:\,\xi_{\s}<\textstyle \frac{r}{t}<\infty\},\] 
where $\xi_{\s}>0$ is yet to be determined. The relevant solution $\tilde U(\xi)$ of 
\eq{U_ode} (i.e., which is defined for $\xi\in (\xi_{\s},\infty)$) 
should give a radial speed $u(t,r)$ which is continuous across $\{t=0\}$.
It follows that $\tilde U(\xi)$ must be the solution to \eq{U_ode} which approaches the 
value $U^*=U_\kaa(-\infty)$ as $\xi\uparrow\infty$. 

Now, as we integrate along decreasing $\xi$-values, in from $\xi=\infty$, the solution
$\tilde U(\xi)$ remains below the solution $-U_\kaa(-\xi)$. This follows since the latter 
function is a solution of \eq{U_ode}
(recall that solutions of \eq{U_ode} lie symmetrically about the origin), and that it 
starts out from $\xi=\infty$ with the value $-U^*>0>U^*=\tilde U(\infty)$.
As a consequence we have that the solution $\tilde U(\xi)$ intersects the 
straight line $l_-$ at some  $\xi$-value $\xi^*$ with $0<\xi^*<-\xi_\w$.

Finally, to determine the shock location $\xi_\s$ we argue as follows. Returning to the solution 
$\hat U(\xi)$ introduced earlier, but now considered for $\xi\in(0,-\xi_\w]$, we let $\hat{\mathcal H}$
denote its associated ``Hugoniot locus.'' That is, $\hat{\mathcal H}$ is the set (curve) of points $(\xi,\hat H(\xi))$
that connect to a point on the solution curve $(\xi,\hat U(\xi))$ through a jump discontinuity
with $U_-=\hat U(\xi)$ and $U_+=\hat H(\xi)$. According to \eq{+ito-}${}_1$,
$\hat{\mathcal H}$ is the graph of the function 
\[\hat H(\xi):=\xi+\frac{a^2}{\hat U(\xi)-\xi}\qquad\text{for $0<\xi<-\xi_\w$.}\]
The following claim follows directly from the properties of the solution $\hat U(\xi)$.
%%%%%%%%%%%%%
\begin{claim}\label{claim}
	The function $\hat H(\xi)$ has the following properties:
	\begin{itemize}
		\item [(i)] $\hat H(\xi)<\xi-a$ for $0<\xi<-\xi_\w$,
		\item [(ii)] $\lim_{\xi\downarrow 0}\hat H(\xi)=-\infty$, and
		\item [(iii)] $\hat H(-\xi_\w)=-U_\w$.
	\end{itemize}
\end{claim}
%%%%%%%%%%%%%
\noindent In particular, it follows from these properties that the graphs of $\hat H(\xi)$ and $\tilde U(\xi)$ 
intersect for some $\xi=\xi_\s\in (0,-\xi_\w)$. (Numerical plots indicate that 
$\hat H(\xi)$ is  strictly increasing on $(0,-\xi_\w)$; if so, $\xi_\s$ is 
uniquely determined. However, we have not been able to provide an analytic proof for this.)
It follows from part (i) of Claim \ref{claim} that the point of intersection lies below $l_-$. Since the graph
of $\hat U(\xi)$ lies between $l_-$ and $l_+$ for $\xi\in (0,-\xi_\w)$, we  
conclude from \eq{isoth_e_1}-\eq{isoth_e_2} that the jump discontinuity with 
with $U_-=\hat U(\xi_\s)$ and $U_+=\hat H(\xi_s)=\tilde U(\xi_\s)$ satisfies the 
entropy condition for a 2-shock. 
See Figure 1.

\bigskip
\noindent
{\bf Summing up:} The radial speed $u(t,r)$ is defined in terms of the  solutions 
$\hat U$, $U_\kaa$, and $\tilde U$ of the similarity ODE \eq{U_ode}, as follows:
\beq\label{u_final}
	u(t,r)=U(\textstyle\frac{r}{t}) :=\left\{\begin{array}{ll}
		\hat U(\frac{r}{t}) & \xi_\w\leq \frac{r}{t}< \xi_\s\\\\
		U_\kaa(\frac{r}{t}) & -\infty<\frac{r}{t}\leq \xi_\w \\\\
		\tilde U(\frac{r}{t}) & \xi_\s<\frac{r}{t}< \infty.
	\end{array}\right.
\eeq
We note that requirement (A) above is met (since $\hat U(0)=0$). Furthermore, this solution 
contains a converging weak discontinuity (``kink'') 
propagating with constant speed along $\{\frac{r}{t}=\xi_\w\}$ for $t<0$ (i.e., $u$ is continuous while 
its first derivatives jump there), and an expanding, entropy admissible 2-shock 
discontinuity propagating with constant speed along $\{\frac{r}{t}=\xi_\s\}$ for $t>0$.
For later use, we record that the radial speed at time of collapse $t=0$ takes 
the constant value
\beq\label{speed_at_collapse}
	u(0,r)\equiv U^*=U_\kaa(-\infty)\qquad\text{for $r>0$.}
\eeq

%%%%%%%%%%%%%%%%%%%%%%%%%%%%%%%%%%%%%%%%%%%
\begin{remark}
	The function $\tilde U(\xi)$ is strictly decreasing on $(\xi_\s,\infty)$ and 
	tends to $U^*<0$ as $\xi\to\infty$. Numerical calculations show that there 
	are cases for which $\tilde U(\xi_\s)>0$ (e.g., this is the case when
	$m=2$, $\beta=-1$), showing that stagnation (vanishing flow velocity) 
	may occur upstream of the expanding shock.
\end{remark}
%%%%%%%%%%%%%%%%%%%%%%%%%%%%%%%%%%%%%%%%%%%

%%%%%%%%%%%%%%%%%%%%%%%%%%%%%%%%%%%%%%%%%%%
\begin{remark}
	In the construction above of $U(\xi)$ on $(-\infty,\xi_\w)$ we made
	use of the particular ``kink'' solution $U_\kaa(\xi)$. We note that, 
	having established that $U_\kaa(-\infty)<0$, we could just as well have
	used any other solution $\check U(\xi)$ of \eq{U_ode} that is located within 
	the region $\mathcal U$ and which exits $\mathcal U$ 
	at a point along the line $\omega$. As noted above, any such 
	solution $\check U(\xi)$ connects smoothly at $\xi=\xi_\w$ to the solution $\hat U(\xi)$ 
	on $[\xi_\w,0]$, and will therefore give converging flows without any weak discontinuities. 
	As $U^*=U_\kaa(-\infty)<0$, it follows that any such solution $\check U(\xi)$ tends to a finite 
	value, $U^{**}$ say, as $\xi\to-\infty$, where $U^{**}<U^*<0$. Then, starting from $U^{**}$ at
	$\xi=+\infty$ and integrating toward the origin, we would generate a 
	solution $U^\circ(\xi)$ (instead of $\tilde U(\xi)$ as above), which
	again could be connected via a jump discontinuity to the solution $\hat U(\xi)$ 
	on $[0,-\xi_\w]$. In particular, we may arrange that $U^{**}$ is so large
	negative that $U^\circ(\xi)$ intersects the Hugoniot curve $\hat H(\xi)$ below 
	the $\xi$-axis; if so, no stagnation occurs in the corresponding flow. 
\end{remark}
%%%%%%%%%%%%%%%%%%%%%%%%%%%%%%%%%%%%%%%%%%%

%%%%%%%%%%%%%%%%%%%%%%%%%%%%%%%%%%%%%%%%%%%%
\subsection{The radial density field $\rho(t,r)$}
%%%%%%%%%%%%%%%%%%%%%%%%%%%%%%%%%%%%%%%%%%%%
With the radial speed defined for all $r\geq 0$ and $t\in\RR$, we turn to the density  
which is given via \eq{fncl_relns}${}_1$,
\beq\label{rho}
	\rho(t,r)=\sgn(t)|t|^\beta\Omega(\xi)\qquad\qquad \xi=\frac{r}{t},
\eeq
where $\Omega$ solves the ODE \eq{isoth_omega_ode} 
\beq\label{Omega_ode_u}
	\frac{\Omega'(\xi)}{\Omega(\xi)}=-\frac{1}{a^2}(U(\xi)-\xi)U'(\xi),
\eeq
and $U(\xi)$ is given by \eq{u_final}. We need to argue that this ODE, together with 
the jump relations at $\xi_\s$, yield a physically
acceptable density field $\rho(t,r)$ satisfying the requirements (B) and (C) in Section 
\ref{conv_div_isothermal}.

As $\beta<0$, it is clear from the second part of requirement (B) that a necessary condition 
on $\Omega$ is that $\Omega(\pm\infty)=0$. However, this is not sufficient to guarantee 
that (B) holds, and we can therefore not use this as an initial condition for the 
$\Omega$-solution. Instead, as we shall see, we can freely assign $\Omega(0-)$ to be any
negative constant $\Omega_0<0$. Having fixed $\Omega_0<0$ we then want to solve the 
ODE \eq{Omega_ode_u}, where $U(\xi)$ is given by \eq{u_final}. 

Before considering the details we outline the order of the various steps for constructing $\Omega(\xi)$.
In what follows, $U(\xi)$ is always given by \eq{u_final}.
We first solve \eq{Omega_ode_u} for $\xi\in [\xi_\w,0]$, obtaining the solution $\hat \Omega(\xi)$ with 
the initial condition $\Omega(0-)=\Omega_0<0$.
We then solve \eq{Omega_ode_u} for $\xi\in (-\infty,\xi_\w]$ with $\Omega(\xi_\w)$ as initial data at $\xi=\xi_\w$, 
obtaining the solution $\Omega_\kaa(\xi)$. As for the velocity $U(\xi)$, the resulting function $\Omega(\xi)$ 
for $\xi\in(-\infty,0)$ suffers a weak discontinuity across $\xi=\xi_\w$. Below we shall show that $\Omega_\kaa(\xi)$ tends to 
zero as $\xi\to-\infty$, and furthermore that it does so in such a manner that
\beq\label{rho_asymp}
	\lim_{t\uparrow0}\rho(t,r)=-C_-r^\beta,
\eeq
where $C_-<0$ is a constant; see \eq{lim_from_below}. 
This will ensure that the constraint (B) is satisfied for times approaching 
zero from below. Since $\beta<0$, it also demonstrates that the density 
field we construct suffers blowup at the origin.

We next need to solve for the density field $\rho(t,r)$ for $\xi_\s<\xi<\infty$, 
and for this it is convenient to switch to the independent variable 
\[x:=\frac{1}{\xi}=\frac{t}{r},\]
and set
\[D(x):=\Omega(\xi).\]
To select the relevant $D$-solution we linearize the ODE for $D(x)$
about the origin in the $(x,D)$-plane and observe that this is a node.
The leading order behavior of the solutions near the origin are of 
the form 
\[D(x)\sim C |x|^{|\beta|},\qquad \text{$C$ constant.}\]
In terms of $\rho(t,r)$ this implies that 
\[\lim_{t\downarrow0}\rho(t,r)=C_+r^\beta,\]
for a constant $C_+$.
Continuity of the density field $\rho(t,r)$ across $\{t=0\}$ requires that we choose $C_+=-C_-$,
where $C_-$ is as in \eq{rho_asymp}. This choice fixes a unique $D$-solution $\tilde D(x)$ for 
$x\gtrsim 0$, which is then unproblematic to extend to all of $[0,x_\s)$, where $x_\s=\frac{1}{\xi_\s}$.
Switching back to  $\xi$ as independent variable, we set
\[\tilde \Omega(\xi):=\tilde D(\textstyle\frac{1}{\xi})\qquad\text{for $\xi_\s<\xi<\infty$.}\]

In particular, this provides us with the value $\tilde\Omega(\xi_\s+)$ at the immediate 
outside of the expanding shock-wave propagating along $\xi=\xi_\s$.
Applying the Rankine-Hugoniot condition \eq{-ito+}${}_2$ with $\bar\xi=\xi_\s$ 
and $\Omega_\pm=\Omega(\xi_\s\pm)$,
we thus determine $\Omega(\xi_\s-)$. This, finally, provides the initial data at $\xi=\xi_\s-$ for 
the relevant solution $\hat\Omega(\xi)$ of \eq{Omega_ode_u} for $\xi\in (0,\xi_\s)$. 
This last step of solving \eq{Omega_ode_u} on $(0,\xi_\s)$ is unproblematic and yields
a final limiting value 
\[\Omega_0'=\lim_{\xi\downarrow 0}\hat\Omega(\xi).\]
We note that differently from the velocity $\hat U(\xi)$, which takes the value zero at $\xi=0$, 
the function $\hat\Omega(\xi)$ will suffer a jump discontinuity there.  
Finally, it is easily verified that the resulting density field satisfies $\rho(t,r)>0$ for all 
$t\in \RR$, $r>0$.
See Figure 2.
We proceed with the details.

%%%%%%%%%%%%%%%%%%
%	FIGURE 
%%%%%%%%%%%%%%%%%%
\begin{figure}\label{Figure_2}
	\centering
	\includegraphics[width=16cm,height=8cm]{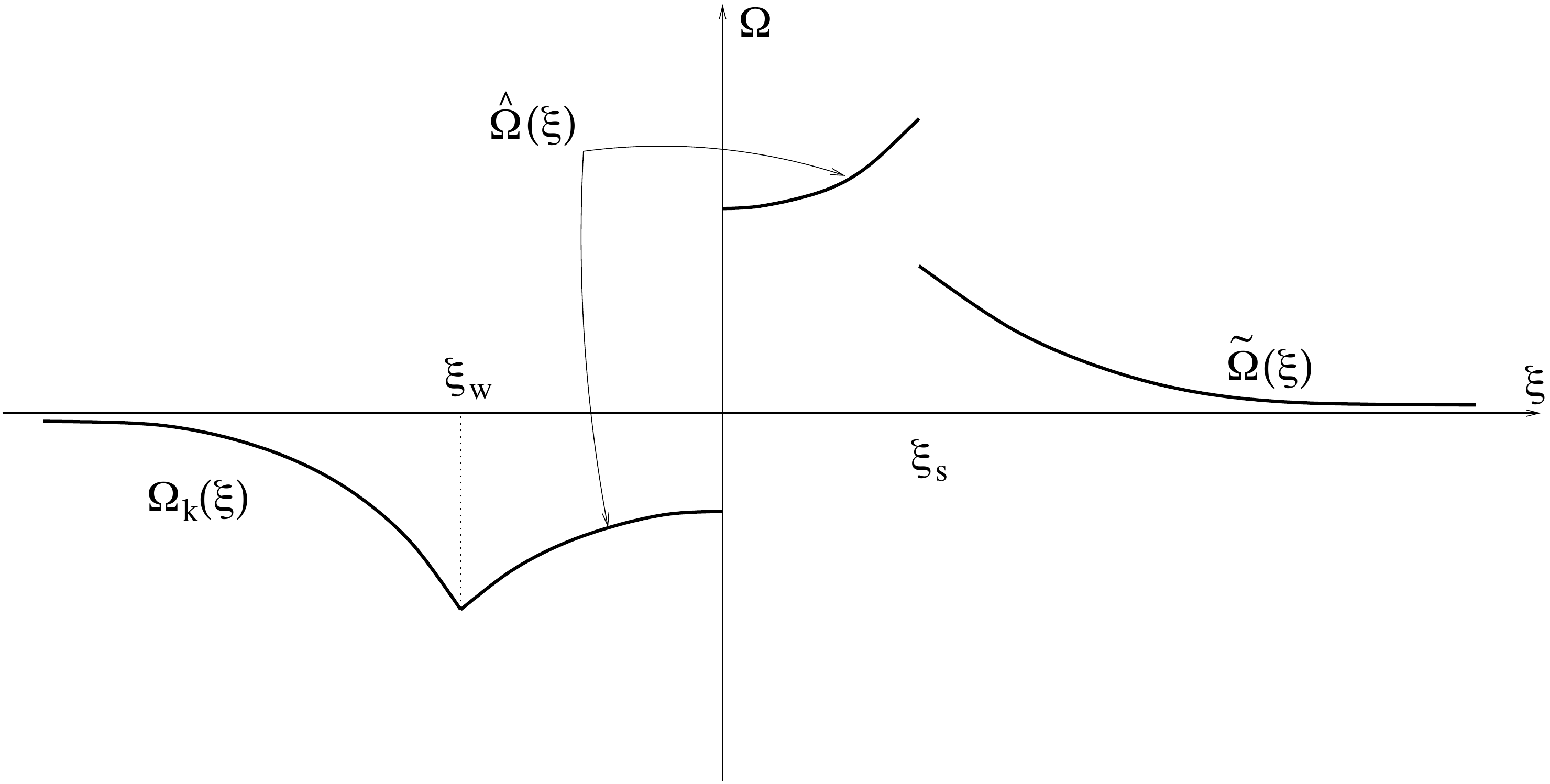}
	\caption{Complete $U(\xi)$-profile (schematic).}
\end{figure}

%%%%%%%%%%%%%%%%%%%%%%%%%%%%%%%%%%%%%%%%%%%%
\subsection{Asymptotics of the density $\rho(t,r)$ for $t\leq 0$}
%%%%%%%%%%%%%%%%%%%%%%%%%%%%%%%%%%%%%%%%%%%%
The first step is to solve 
\beq\label{hat_Omega_ode_1}
	\frac{\hat\Omega'(\xi)}{\hat\Omega(\xi)}=-\frac{1}{a^2}(\hat U(\xi)-\xi)\hat U'(\xi)=:\hat F(\xi)
	\qquad\text{for $\xi\in[\xi_\w,0]$,}
\eeq
where $\hat U(\xi)$ was determined above. As initial data we fix any constant $\Omega_0<0$ 
and set 
\[\hat \Omega(0-):=\Omega_0.\]
It follows from the properties of $\hat U(\xi)$ that $\hat F(\xi)$ is a bounded, smooth 
function on $[\xi_\w,0]$, such that solving \eq{hat_Omega_ode_1} is  unproblematic. We 
note that 
\beq\label{Omega_sign1}
	\hat \Omega(\xi)< 0\qquad\text{and}\qquad \hat \Omega'(\xi)\geq 0\qquad\text{for $\xi\in[\xi_\w,0]$.} 
\eeq
Next we want to solve 
\beq\label{Omega_k_ode}
	\frac{\Omega_\kaa'(\xi)}{\Omega_\kaa(\xi)}=-\frac{1}{a^2}(U_\kaa(\xi)-\xi) U'_\kaa(\xi)=: F_\kaa(\xi)
	\qquad\text{for $\xi\in(-\infty,\xi_\w]$,}
\eeq
where $U_\kaa(\xi)$ was determined above. To establish \eq{rho_asymp} we first show that 
\beq\label{est_0}
	\int_{-\infty}^{\xi_\w} |F_\kaa(\eta)-\textstyle\frac{\beta}{\eta}|\,d\eta<\infty.
\eeq
Indeed, by using that
\[U'_\kaa=\frac{a^2}{(U_\kaa-\xi)^2-a^2}\Big(\beta+\frac{mU_\kaa}{\xi}\Big),\]
together with the fact that $U_\kaa(\xi)\to U^*=U_\kaa(-\infty)<0$, it is straightforward to verify that 
\beq\label{est_1}
	|F_\kaa(\xi)-\textstyle\frac{\beta}{\xi}|\leq\frac{C}{\xi^2}\qquad\text{for $\xi\in(-\infty,\xi_\w]$,}
\eeq
for a suitable constant $C$, and \eq{est_0} follows. Integrating \eq{Omega_k_ode}, we obtain
\[\Omega_\kaa(\xi)=\Omega_\w\frac{|\xi|^\beta}{|\xi_\w|^\beta}
\cdot\exp\Big(\int_{\xi}^{\xi_\w} \textstyle\frac{\beta}{\eta}-F_\kaa(\eta)\,d\eta\Big),\]
where $\Omega_\w:=\Omega_\kaa(\xi_\w)<0$. Applying \eq{est_1} yields
\beq\label{Omega_asymp}
	\Omega_\kaa(\xi) \sim C_-|\xi|^\beta\qquad\text{as $\xi\to-\infty$,}
\eeq
where 
\[C_-=\frac{\Omega_\w}{|\xi_\w|^\beta}
\cdot\exp\Big(\int_{\xi}^{\xi_\w} \textstyle\frac{\beta}{\eta}-F_\kaa(\eta)\,d\eta\Big)<0.\]
Applying this in \eq{rho} we obtain
\beq\label{lim_from_below}
	\lim_{t\uparrow0}\rho(t,r)=\lim_{t\uparrow0}\,\,\sgn(t)|t|^\beta\Omega_\kaa(\textstyle\frac{r}{t})
	=-C_-r^\beta\qquad \text{at any fixed location $r>0$,}
\eeq
verifying \eq{rho_asymp}. We also note that \eq{Omega_k_ode}, together with 
the properties of $U_\kaa(\xi)$, imply that
\beq\label{Omega_sign2}
	\Omega_\kaa(\xi)<0\qquad\text{and}\qquad  \Omega_\kaa'(\xi)<0\qquad\text{for $(-\infty,\xi_\w)$.}
\eeq

%%%%%%%%%%%%%%%%%%%%%%%%%%%%%%%%%%%%%%%%%%%%
\subsection{The density $\rho(t,r)$ for $t\geq 0$.}
%\label{rho_for_pos_t}
%%%%%%%%%%%%%%%%%%%%%%%%%%%%%%%%%%%%%%%%%%%%
To identify the relevant solution $\tilde \Omega(\xi)$ for $\xi\in (\xi_\s,\infty)$,
we switch to the independent variable $x=\frac{1}{\xi}$ and set $\tilde D(x)=\tilde \Omega(\frac{1}{x})$.
The ODE for $\tilde D(x)$ is given by \eq{Omega_ode_u} 
\beq\label{D_ode}
	\frac{\tilde D'(x)}{\tilde D(x)}=\frac{1}{a^2x^2}
	\big[\textstyle \tilde U\big(\frac{1}{x}\big)-\frac{1}{x}\big]\tilde U'\big(\frac{1}{x}\big)
	\qquad \text{for $0<x<x_\s,$}
\eeq
where $\tilde U(\xi)$ was determined above. It follows from requirement (B) in Section 
\ref{conv_div_isothermal} that we must have 
$\tilde D(0)=0$. Linearizing \eq{D_ode} about $(x,\tilde D)=(0,0)$ shows that
the origin is a node where 
\beq\label{asymp_2}
	\tilde D(x)\sim C_+ x^{-\beta}\qquad \text{for $x\gtrsim 0$,}
\eeq
or 
\[\tilde \Omega(\xi)\sim C_+ \xi^{\beta}\qquad \text{as $\xi\to+\infty$.}\]
This gives 
\beq\label{lim_from_above}
	\lim_{t\downarrow 0} \rho(t,r)=\lim_{t\downarrow0}\,\,t^\beta\tilde \Omega(\textstyle\frac{r}{t})
	=C_+r^\beta\qquad \text{at any fixed location $r>0$.}
\eeq
%[Note: the same linearization shows that \eq{asymp_2}, with $x$ replaced by $|x|$, remains to hold
%for $x\lesssim0$. However, in order to use this we would still have to argue that the relevant solution 
%$D_\kaa(x)=\Omega_\kaa(\frac{1}{x})$ of \eq{D_ode} approaches zero as $x\uparrow 0$. It appears 
%easier to argue for this in terms of $\Omega_\kaa(\xi)$ as we did above.]
Comparing with \eq{lim_from_below} and imposing continuity of $\rho(t,r)$ across $\{t=0\}$,
implies that $C_+=-C_-$, and this selects the unique, relevant solution $\tilde D(x)$ for $x\gtrsim 0$.

It is now unproblematic to integrate \eq{D_ode} for $x\in (0,x_\s)$ (where $x_\s=\frac{1}{\xi_\s}$), and
it follows from \eq{D_ode}, together with the properties of $\tilde U(\xi)$, \eq{asymp_2}, and $C_+>0$, that 
$\tilde D(x)>0$ and $\tilde D'(x)>0$
for $0<x<x_\s$. We therefore obtain that
\beq\label{Omega_sign3}
	\tilde \Omega(\xi)>0\qquad\text{and}\qquad \tilde \Omega'(\xi)<0\qquad\text{for $\xi\in(\xi_\s,\infty)$.}
\eeq

Having obtained $\tilde \Omega(\xi)$ for $\xi>\xi_\s$, we use the Rankine-Hugoniot relation 
\eq{-ito+}${}_2$ with $\bar \xi=\xi_\s$, $\Omega_+=\tilde \Omega(\xi_\s)$ and $U_+=\tilde U(\xi_\s)$, 
to calculate $\Omega_-$. This last value is used as initial data at $\xi=\xi_\s$ for the ODE 
\beq\label{hat_Omega_ode_2}
	\frac{\hat\Omega'(\xi)}{\hat\Omega(\xi)}=-\frac{1}{a^2}(\hat U(\xi)-\xi)\hat U'(\xi)
	\qquad\text{for $\xi\in(0,\xi_\s)$.}
\eeq
We note that, since $(\tilde U(\xi_\s)-\xi_\s)^2>a^2$, \eq{-ito+}${}_2$ gives
\[\hat \Omega(\xi_\s)>\tilde \Omega(\xi_\s)>0.\]
It then follows from the properties of $\hat U(\xi)$ that the right-hand side of 
\eq{hat_Omega_ode_2} is a bounded and positive function on $[0,\xi_\s]$. 
Consequently, $\hat\Omega(\xi)$ is increasing there and approaches a 
strictly positive value $\Omega_0'$ at $\xi=0+$:
\beq\label{Omega_sign4}
	\hat \Omega(\xi)>0\qquad\text{and}\qquad \hat\Omega'(\xi)>0
	\qquad\text{for $\xi\in(0,\xi_\s)$,}
	\qquad\text{and}\qquad\lim_{\xi\downarrow 0}\hat \Omega(\xi)=\Omega_0'>0.
\eeq

\bigskip
\noindent
{\bf Summing up:} The density field $\rho(t,r)$ is defined in terms of the  solutions 
$\hat \Omega$, $\Omega_\kaa$, and $\tilde \Omega$ of the similarity ODE \eq{Omega_ode} 
as determined above, as follows:
\beq\label{rho_final}
	\rho(t,r)=\sgn(t)|t|^\beta\Omega(\textstyle\frac{r}{t}) :=\left\{\begin{array}{ll}
		- |t|^\beta\hat\Omega(\frac{r}{t}) & \xi_\w\leq \frac{r}{t}\leq 0\\\\
		- |t|^\beta \Omega_\kaa(\frac{r}{t}) & -\infty<\frac{r}{t}\leq \xi_\w \\\\
		t^\beta\tilde\Omega(\frac{r}{t}) & \xi_\s<\frac{r}{t}< \infty \\\\
		t^\beta\hat \Omega(\frac{r}{t}) &0\leq \frac{r}{t}< \xi_\s.\\
	\end{array}\right.
\eeq
We note that, as for the radial speed given by \eq{u_final},  the density field suffers a weak discontinuity
across $\{\frac{r}{t}=\xi_\w\}$ for $t<0$, and a jump discontinuity across $\{\frac{r}{t}=\xi_\s\}$ for $t>0$.
As detailed at the end of Section \ref{u_for_pos_t}, the resulting shock wave along $\{r=\xi_\s t\}$ 
is, by construction, an entropy admissible 2-shock for the isothermal Euler system.
Next, recalling \eq{Omega_sign1}, \eq{Omega_sign2}, \eq{Omega_sign3}, and \eq{Omega_sign4}, 
we have that $\Omega(\xi)\neq 0$ for all values of $\xi$. Furthermore, the density field at the time 
of collapse $t=0$ is given by
\beq\label{density_at_collapse}
	\rho(0,r)=|C_-|r^\beta\qquad r>0.
\eeq
It follows from this that requirement (C) above is met by the density field given by
\eq{rho_final}: $\rho(t,r)>0$ for all $t\in\RR$ and all $r\geq 0$.
Finally, \eq{lim_from_below}, \eq{lim_from_above}, and the choice $C_+=-C_-$, show 
that also requirement (B) is satisfied.

%%%%%%%%%%%%%%%%%%%%%
\begin{remark}
	The above construction of $\rho(t,r)$ and $u(t,r)$ provides a 2-parameter family of 
	concrete solutions to the radial, isothermal Euler system in $n=2$ and $n=3$ space 
	dimensions. The solutions depend on the similarity exponent $\beta$, which varies 
	in $(-n+1,0)$ so as to satisfy Assumption \ref{assmpn}, and on the constant $\Omega_0<0$, 
	which determines the density along the center of motion before collapse 
	($\rho(t,0)=|\Omega_0||t|^\beta$ for $t<0$).
\end{remark}
%%%%%%%%%%%%%%%%%%%%%

%%%%%%%%%%%%%%%%%%%%%%%%%%%%%%%%%%%%%%%%%%%%%
\section{Weak and radial weak Euler solutions}
\label{weak_solns} 
%%%%%%%%%%%%%%%%%%%%%%%%%%%%%%%%%%%%%%%%%%%%% 
It remains to verify that the radial solutions of the isothermal Euler system 
constructed above do indeed provide genuine, weak solutions
to the original, multi-d isothermal Euler system 
\eq{mass_m_d_isthrml_eul}-\eq{mom_m_d_isthrml_eul}. 
In this section we formulate the definition of a weak solution to the 
barotropic Euler system: first for general, multi-d solutions, and 
then specialized to the case of radial solutions.

%%%%%%%%%%%%%%%%%%%%%%%%%%%%%%%%%%%%%%%%%%%%%
\subsection{Multi-d weak solutions}
%\label{multi-d_weak_solns} 
%%%%%%%%%%%%%%%%%%%%%%%%%%%%%%%%%%%%%%%%%%%%%
We write $\rho(t)$ for $\rho(t,\cdot)$ etc., ${\bf u}=(u_1,\dots,u_n)$, 
$u:=|{\bf u}|$, and let ${\bf x}=(x_1,\dots,x_n)$ denote the spatial variable in $\RR^n$,
while $r=|{\bf x}|$ varies over $\RR_0^+=[0,\infty)$.

%%%%%%%%%%%%%%%%%%%%%%%%%%%%
\begin{definition}\label{weak_soln}
	Consider the compressible, isothermal Euler system 
	\eq{mass_m_d_isthrml_eul}-\eq{mom_m_d_isthrml_eul}
	in $n$ space dimensions with a given pressure function $p=p(\rho)\geq 0$.
	Then the measurable functions $\rho,\, u_1,\dots,u_n:\RR_t\times \RR_{\bf x}^n\to \RR$ 
	constitute a {\em weak solution} to \eq{mass_m_d_isthrml_eul}-\eq{mom_m_d_isthrml_eul} 
	provided that:
	\begin{itemize}
		\item[(1)] the maps $t\mapsto \rho(t)$ and $t\mapsto \rho(t) u(t)$
		belong to $C^0(\RR_t;L^1_{loc}(\RR^n_{\bf x}))$;
		\item[(2)]  the functions $\rho u^2$ and $p$
		belong to $L^1_{loc}(\RR_t\times\RR^n_{\bf x})$;
		\item[(3)] the conservation laws for mass and momentum are 
		satisfied weakly in sense that
		\beq\label{m_d_mass_weak}
			\int_\RR\int_{\RR^n} \rho\vp_t+\rho{\bf u}\cdot\nabla_{\bf x}\vp
			\, d{\bf x}dt =0
		\eeq
		and
		\beq\label{m_d_mom_weak}
			\int_\RR\int_{\RR^n} \rho u_i\vp_t
			 +\rho u_i{\bf u}\cdot\nabla_{\bf x}\vp+p\vp_{x_i}\, d{\bf x}dt =0
			 \qquad \text{for $i=1,\dots, n$,} 
		\eeq
		whenever $\vp\in C_c^1(\RR_t\times \RR^n_{\bf x})$ ($C^1$ functions with compact support).
	\end{itemize}
\end{definition}
%%%%%%%%%%%%%%%%%%%%%%%%%%%%%%%%
\begin{remark}
	Here, condition (1) guarantees that the conserved quantities define
	continuous maps into $L^1_{loc}(\RR^n_{\bf x})$, which is the natural function
	space in this setting. Taken together, conditions (1) and (2) ensure that all terms occurring 
	in the weak formulations \eq{m_d_mass_weak} and \eq{m_d_mom_weak} are 
	locally integrable in space and time.
\end{remark}
%%%%%%%%%%%%%%%%%%%%%%%%%%%%%%%%
%%%%%%%%%%%%%%%%%%%%%%%%%%%%%%%%
\begin{remark}
	Our goal is to show that the converging-diverging flow 
	\beq\label{assmbld}
		\rho(t,{\bf x})=\rho(t,r),\qquad {\bf u}(t,{\bf x})=u(t,r)\frac{{\bf x}}{r}, 
	\eeq
	where $\rho(t,r)$ and $u(t,r)$ are given by
	\eq{rho_final} and \eq{u_final}, respectively, constitute a weak solution to 
	\eq{mass_m_d_isthrml_eul}-\eq{mom_m_d_isthrml_eul} (with $p=a^2\rho$)
	according to the definition above. Since these flows by construction involve a single, 
	compressive shock wave, we do not address admissibility of weak solutions.
\end{remark}
%%%%%%%%%%%%%%%%%%%%%%%%%%%%%%%%

%%%%%%%%%%%%%%%%%%%%%%%%%%%%%%%
\subsection{Radial weak solutions}
%\label{rad_weak_solns}
%%%%%%%%%%%%%%%%%%%%%%%%%%%%%%%
We next rewrite Definition \ref{weak_soln} for radial solutions.
For this we use the following notation. As above $m:=n-1$ 
and we set 
\[\RR^+=(0,\infty),\qquad \RR_0^+=[0,\infty),\qquad 
L^1_{(loc)}(dt\times r^mdr)=L^1_{(loc)}(\RR\times\RR^+_0,dt\times r^mdr).\]
Also, $C^1_c(\RR\times\RR^+_0)$ denotes the set of real-valued functions 
$\psi(t,r)$ defined on $\RR\times\RR^+_0$ and with the property that $\psi$ 
is $C^1$ smooth on $\RR\times\RR^+_0$ and vanishes outside
$[-\bar t,\bar t]\times[0,\bar r]$ for some $\bar t,\, \bar r\in\RR^+$.
Finally, we let $C^1_0(\RR\times\RR^+_0)$ denote the set of those functions 
$\psi\in C^1_c(\RR\times\RR^+_0)$ with the additional property that 
$\psi(t,0)\equiv 0$.

Using these function classes, the weak formulation of the  multi-d 
Euler system \eq{mass_m_d_isthrml_eul}-\eq{mom_m_d_isthrml_eul},
for radial solutions, takes the following form.

%%%%%%%%%%%%%%%%%%%%%%%%%%%%
\begin{definition}\label{rad_symm_weak_soln}
	Consider the radial version \eq{mass}-\eq{momentum} of the compressible Euler 
	system  \eq{mass_m_d_isthrml_eul}-\eq{mom_m_d_isthrml_eul} 
	with a given pressure function $p=p(\rho)\geq 0$.
	Then the measurable functions $\rho,\, u:\RR_t\times \RR^+_r\to \RR$  
	constitute a {\em radial weak solution} to \eq{mass}-\eq{momentum} provided that:
	\begin{itemize}
		\item[(i)] the maps $t\mapsto \rho(t)$ and $t\mapsto \rho(t)u(t)$
		belong to $C^0(\RR_t;L^1_{loc}(r^mdr))$;
		\item[(ii)]  the functions $\rho u^2$ and $p$ belong to 
		$L^1_{loc}(dt\times r^mdr)$;
		\item[(iii)] the conservation laws for mass and momentum are 
		satisfied in the sense that
		\begin{align}
			\int_{\RR}\int_{\RR^+} \left(\rho\psi_t+\rho u\psi_r\right) r^mdrdt &=0 
			\qquad\forall \psi\in C^1_c(\RR\times\RR^+_0) \label{radial_mass_weak}\\
			\int_{\RR}\int_{\RR^+} \left(\rho u\psi_t
			+\rho u^2\psi_r+p\big(\psi_r+\textstyle\frac{m\psi}{r}\big)\right) r^mdrdt &=0 
			\qquad\forall \psi\in C^1_0(\RR\times\RR^+_0).\label{radial_mom_weak}
		\end{align}
	\end{itemize}
\end{definition}
%%%%%%%%%%%%%%%%%%%%%%%%%%%%
The demonstration that a radial weak solution $(\rho,u)$ yields, 
via \eq{assmbld}, a weak solution of the multi-d system according to 
Definition \ref{weak_soln}, was provided by Hoff  \cite{hoff} in the context 
of radial, isentropic Navier-Stokes flows. (See \cite{jt1} for the corresponding 
analysis in the case of radial, non-isentropic Euler flows).

%%%%%%%%%%%%%%%%%%%%%%%%%%%%%%%%%%%%%%
%%%%%%%%%%%%%%%%%%%%%%%%%%%%%%%%%%%%%%
\section{Radial converging-diverging similarity solutions as  weak solutions}
\label{sim_weak_solns} 
%%%%%%%%%%%%%%%%%%%%%%%%%%%%%%%%%%%%%
%%%%%%%%%%%%%%%%%%%%%%%%%%%%%%%%%%%%%%
In this section we return to isothermal flow ($p=a^2\rho$) and the radial 
converging-diverging 
similarity solutions constructed in Section \ref{constr_conv_div_solns}.
We want to establish properties (i), (ii), and (iii) in Definition 
\ref{rad_symm_weak_soln} for these solutions, and we first consider the continuity and 
integrability requirements in (i) and (ii). The weak forms of the equations
are treated in Section \ref{weak_forms}. 

%%%%%%%%%%%%%%%%%%%%%%%%%%%%%%%%%%%
\subsection{Continuity and local integrability}
%\label{cont_and_loc_integr}
%%%%%%%%%%%%%%%%%%%%%%%%%%%%%%%%%%%
With $\rho(t,r)$ and $u(t,r)$ given by \eq{rho_final} and \eq{u_final}, we  
proceed to verify parts (i) and (ii) of Definition \ref{rad_symm_weak_soln}. 
For this we fix $\br>0$, define
\[M(t;\bar r):=\int_0^{\br} \rho(t,r)r^m\, dr,\qquad 
I_q(t;\bar r):=\int_0^{\br} \rho(t,r)|u(t,r)|^qr^m\, dr\qquad (q=1,\, 2),\]
and observe that, in the particular case under consideration, where 
$p\propto \rho$, (i) and (ii) both follow once we verify 
that the maps $t\mapsto M(t;\bar r)$, $t\mapsto I_1(t;\bar r)$, and
$t\mapsto I_2(t;\bar r)$ are continuous at all times $t\in\RR$. 
Now, as $\rho(t,r)$ and $u(t,r)$ are bounded functions, except at 
the time of collapse ($t=0$), it is sufficient to verify 
the continuity of $M(t;\bar r)$ and $I_q(t;\bar r)$ ($q=1$, $2$) across $t=0$.

According to \eq{density_at_collapse}, together with the standing 
assumption $\beta+m>0$, we have that $M(0;\br)$ is finite and given by
\[M(0;\br)=\frac{|C_-|}{\beta+n}\bar r^{\beta+n}.\]
For $t<0$ (and small enough that $\xi_\w t<\bar r$) we have 
\begin{align*}
	M(t;\bar r)&=\int_0^{\bar r}\rho(t,r) r^m\, dr
	=\int_0^{\br/t}\sgn(t)|t|^\beta |\Omega(\xi)|(t\xi)^mt\, d\xi\nn\\
	&=|t|^{\beta+n}\Big[\int_{\br/t}^{\xi_\w}|\Omega_\kaa(\xi)||\xi|^m\, d\xi 
	+\int_{\xi_\w}^0|\hat \Omega(\xi)||\xi|^m\, d\xi \Big].
\end{align*}
Here the last term in the brackets is a bounded number, while L'H\^ opital's rule 
applied to the first term gives
\begin{align*}
	\lim_{t\uparrow0}M(t;\bar r)&=\lim_{t\uparrow0}\frac{1}{|t|^{-\beta-n}}
	\int_{\br/t}^{\xi_\w}|\Omega_\kaa(\xi)||\xi|^m\, d\xi\nn\\
	&=\lim_{t\uparrow0}\frac{\br^n}{\beta+n}|t|^\beta|\Omega_\kaa({\textstyle\frac{\br}{t}})|
	=\frac{|C_-|\br^{\beta+n}}{\beta+n},
\end{align*}
where we have used \eq{lim_from_below}. An entirely similar calculation, now using 
\eq{lim_from_above} and with $\xi_\s$ playing the role of $\xi_\w$, shows that 
\[\lim_{t\downarrow0}M(t;\bar r)=\lim_{t\downarrow0}
\frac{\br^n}{\beta+n}t^\beta\tilde\Omega({\textstyle\frac{\br}{t}})
=\frac{C_+\br^{\beta+n}}{\beta+n}.\]
As $C_+=|C_-|$, this establishes the continuity of $t\mapsto M(t;\bar r)$ at time 
$t=0$, and thus for all times.

Next, according to \eq{speed_at_collapse} and \eq{density_at_collapse}, we have
\[I_q(0;\bar r)=\int_0^{\br} \rho(0,r)|u(0,r)|^qr^m\, dr
=\frac{|C_-||U^*|^q}{\beta+n}\br^{\beta+n}.\]
As above, for $t\lesssim 0$, we have
\begin{align*}
	I_q(t;\bar r)&=\int_0^{\bar r}\rho(t,r) |u(t,r)|^q r^m\, dr
	=\int_0^{\br/t}\sgn(t)|t|^\beta |\Omega(\xi)||U(\xi)|^q(t\xi)^mt\, d\xi\nn\\
	&=|t|^{\beta+n}\Big[\int_{\br/t}^{\xi_\w}|\Omega_\kaa(\xi)||U_\kaa(\xi)|^q|\xi|^m\, d\xi 
	+\int_{\xi_\w}^0|\hat \Omega(\xi)||\hat U(\xi)|^q|\xi|^m\, d\xi \Big].
\end{align*}
Again, here the last term in the brackets is a bounded number, while L'H\^ opital's rule 
applied to the first term gives
\begin{align*}
	\lim_{t\uparrow0}I_q(t;\bar r)&=\lim_{t\uparrow0}\frac{1}{|t|^{-\beta-n}}
	\int_{\br/t}^{\xi_\w}|\Omega_\kaa(\xi)||U_\kaa(\xi)|^q|\xi|^m\, d\xi\nn\\
	&=\lim_{t\uparrow0}\frac{\br^n}{\beta+n}|t|^\beta|\Omega_\kaa({\textstyle\frac{\br}{t}})|
	|U_\kaa({\textstyle\frac{\br}{t}})|^q
	=\frac{|C_-||U^*|^q}{\beta+n}\br^{\beta+n},
\end{align*}
where we have used \eq{lim_from_below}. A similar calculation shows that 
\[\lim_{t\uparrow0}I_q(t;\bar r)
	=\frac{C_+|U^*|^q}{\beta+n}\br^{\beta+n},\]
As $C_+=|C_-|$, this establishes the continuity of the maps $t\mapsto I_q(t;\bar r)$,
$q=1$, $2$, at time $t=0$, and thus for all times.

We have thus verified requirements (i) and (ii) of Definition \ref{rad_symm_weak_soln}
for the isothermal converging-diverging solutions $(\rho(t,r),u(t,r))$ constructed in Section  
\ref{constr_conv_div_solns}.

%%%%%%%%%%%%%%%%%%%%%%%%%%%%%%%%%%%
\subsection{Weak form of the equations}\label{weak_forms}
%%%%%%%%%%%%%%%%%%%%%%%%%%%%%%%%%%%
Finally, for part (iii) of Definition \ref{rad_symm_weak_soln}, we need to verify the
weak forms \eq{radial_mass_weak}, \eq{radial_mom_weak}. For this 
we shall exploit that the local integrability properties in parts (i) and (i) of Definition 
\ref{rad_symm_weak_soln} have been verified. 
The issue will then reduce to estimating the fluxes of the conserved quantities
across spheres of vanishing radii.

For $\psi\in C^1_c(\RR\times\RR^+_0)$, with $\supp\psi\subset[-T,T]\times [0,\br]$,
and any small $\delta>0$, we define the regions
\[J_\delta=\left\{(t,r)\,|\, -T<t<T,\, \delta<r<\br,\, \textstyle\frac{t}{r}<\frac{1}{\xi_\s} \right\},\]
and
\[K_\delta=\left\{(t,r)\,|\, -T<t<T,\, \delta<r<\br,\, \textstyle\frac{t}{r}>\frac{1}{\xi_\s} \right\},\]
(see Figure 3), and set
\begin{align}
	M(\psi)&:=\iint_{\RR\times\RR^+} \left(\rho\psi_t+\rho u\psi_r\right) r^mdrdt \nn\\
	&= \Big\{\iint_{\RR\times [0,\delta]} +\iint_{J_\delta} +\iint_{K_\delta}\Big\}
	\left(\rho\psi_t+\rho u\psi_r\right) r^mdrdt \nn\\
	&=:M_\delta(\psi)
	+\Big\{\iint_{J_\delta} +\iint_{K_\delta}\Big\}
	\left(\rho\psi_t+\rho u\psi_r\right) r^mdrdt
	\label{M_psi}
\end{align}
and 
\begin{align}
	I(\psi)&:=\iint_{\RR\times\RR^+}\left(\rho u\psi_t
	+\rho u^2\psi_r+p\big(\psi_r+\textstyle\frac{m\psi}{r}\big)\right) r^mdrdt \nn\\
	&= \Big\{\iint_{\RR\times [0,\delta]} +\iint_{J_\delta} +\iint_{K_\delta}\Big\}
	\left(\rho u\psi_t
	+\rho u^2\psi_r+p\big(\psi_r+\textstyle\frac{m\psi}{r}\big)\right) r^mdrdt \nn\\
	&=:I_\delta(\psi)
	+\Big\{\iint_{J_\delta}+\iint_{K_\delta}\Big\}
	\left(\rho u\psi_t
	+\rho u^2\psi_r+p\big(\psi_r+\textstyle\frac{m\psi}{r}\big)\right) r^mdrdt.
	\label{I_psi}
\end{align}
The goal is to verify that $M(\psi)$ and $I(\psi)$ vanish by showing that the right hand 
sides of \eq{M_psi} and \eq{I_psi} vanish as $\delta\downarrow 0$.

%%%%%%%%%%%%%%%%%%
%	FIGURE 
%%%%%%%%%%%%%%%%%%
\begin{figure}\label{Figure_3}
	\centering
	\includegraphics[width=8cm,height=9cm]{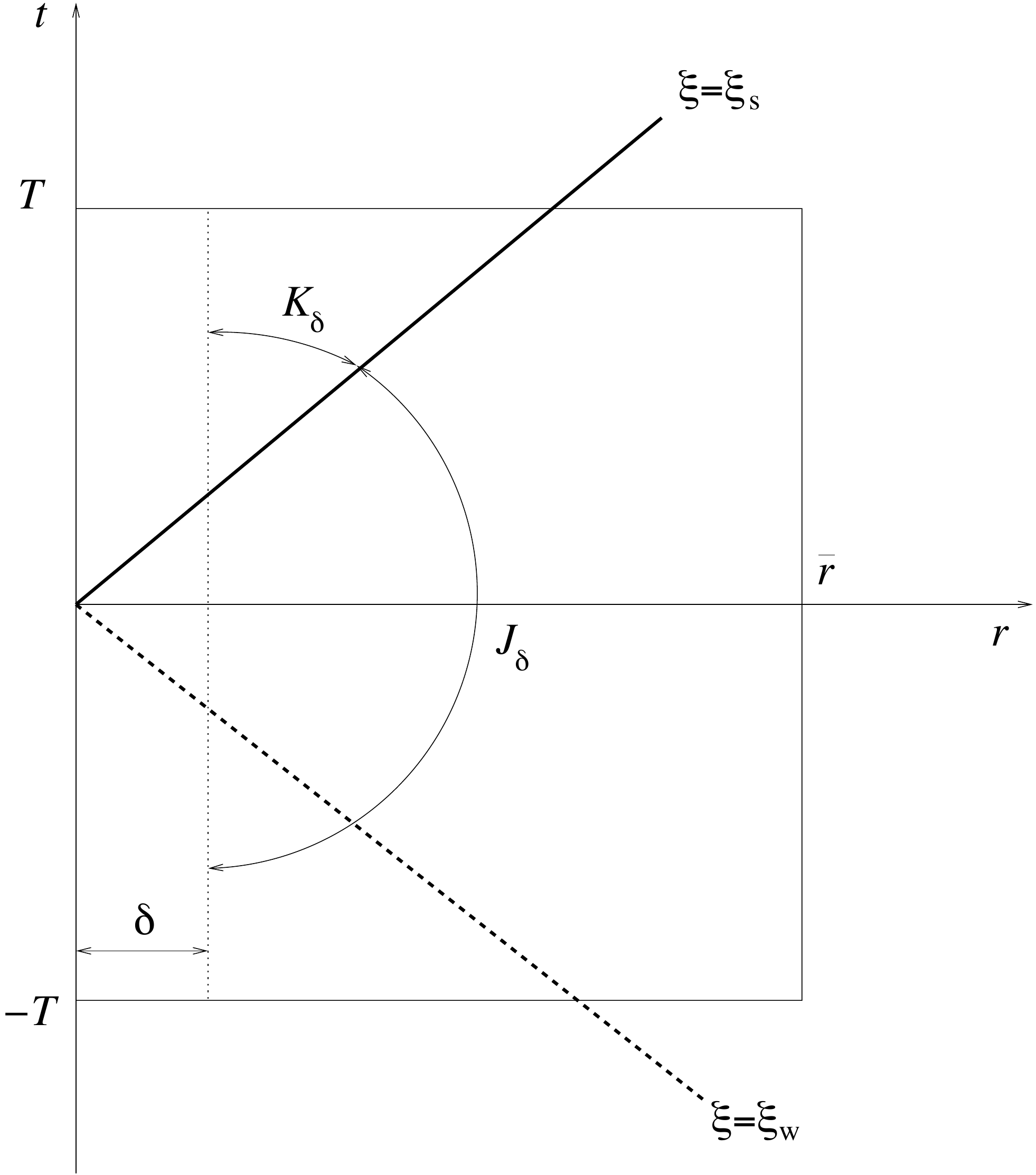}
	\caption{Regions of integration in the weak formulation.}
\end{figure}

We first note that the continuity of the maps $t\mapsto M(t;\bar r)$, 
$t\mapsto I_1(t;\bar r)$, and $t\mapsto I_2(t;\bar r)$, which was established above,
implies the local $r^mdrdt$-integrability of $\rho$, $p\propto \rho$, $\rho u$, and
$\rho u^2$. As a consequence, both $M_\delta(\psi)$ and $I_\delta(\psi)$ tend to zero
as $\delta\downarrow 0$. (Note that for $I_\delta(\psi)$, we make use of the fact that
$\psi$ belongs to the space $C^1_0(\RR\times\RR^+_0)$; in particular,
$\frac{m\psi}{r}$ is a bounded term.)

It remains to estimate the integrals over $J_\delta$ and $K_\delta$ in \eq{M_psi} 
and \eq{I_psi}. For this we first recall that $(\rho,u)$, by construction, is a 
classical (Lipschitz) solution of the isentropic Euler system \eq{mass}-\eq{momentum}
within each of $J_\delta$ and $K_\delta$, and that the Rankine-Hugoniot relations 
\eq{rh_1}-\eq{rh_2}, with $\dot{\mathcal R}=\xi_\s$,
are satisfied across their common boundary along the straight line $\{r=\xi_\s t\}$. 
Applying the divergence theorem to each region we therefore have,
\beq\label{part1}
	\Big\{\iint_{J_\delta} +\iint_{K_\delta}\Big\}\left(\rho\psi_t+\rho u\psi_r\right)\, r^mdrdt
	=-\delta^m\int_{-T}^T(\rho u\psi)(t,\delta)\, dt
\eeq
and
\beq\label{part2}
	\Big\{\iint_{J_\delta}+\iint_{K_\delta}\Big\}\left(\rho u\psi_t
	+\rho u^2\psi_r+p\big(\psi_r+\textstyle\frac{m\psi}{r}\big)\right) r^mdrdt
	=-\delta^m\int_{-T}^T[(\rho u^2+p)\psi](t,\delta)\, dt.
\eeq
Since the speed $u(t,r)$ under consideration is globally bounded,  $\psi(t,r)$ is
a bounded function, and $p\propto \rho$, it follows that to estimate these expressions,
it suffices to consider the single quantity $\delta^m\int_{-T}^T \rho(t,\delta)$. We have,
using \eq{rho_final} and switching to $\xi$ as integration variable,
\begin{align}
	\delta^m\int_{-T}^T \rho(t,\delta)\, dt 
	&= \delta^{n+\beta}\Big\{\int_{\xi_\w}^{-\delta/T}+\int_{-\infty}^{\xi_\w}
	+\int_{\xi_\s}^\infty+\int_{\delta/T}^{\xi_\s}\Big\}\frac{|\Omega(\xi)|}{|\xi|^{\beta+2}}\,d\xi\nn\\
	&= \delta^{n+\beta}\Big\{\int_{\xi_\w}^{-\delta/T}\frac{|\hat\Omega(\xi)|}{|\xi|^{\beta+2}}\,d\xi
	+\int_{-\infty}^{\xi_\w}\frac{|\Omega_\kaa(\xi)|}{|\xi|^{\beta+2}}\,d\xi
	+\int_{\xi_\s}^\infty\frac{\tilde \Omega(\xi)}{\xi^{\beta+2}}\,d\xi
	+\int_{\delta/T}^{\xi_\s}\frac{\hat \Omega(\xi)}{\xi^{\beta+2}}\,d\xi\Big\}.
	\label{flux_at_delta}
\end{align}
According to \eq{Omega_asymp} and \eq{asymp_2}, we have, for a suitable constant $C$,
\[|\Omega_\kaa(\xi)|\leq C|\xi|^\beta \quad\text{for $\xi<\xi_\w$, and}\quad
\tilde\Omega(\xi)\leq C\xi^\beta \quad\text{for $\xi>\xi_\s$.}\]
Also, according to the construction in Section \ref{constr_conv_div_solns}, $\hat\Omega(\xi)$
is a bounded function. Using these in \eq{flux_at_delta}, we get that
\begin{align*}
	\delta^m\int_{-T}^T \rho(t,\delta)\, dt 
	&\leq const. \delta^{n+\beta}
	\left\{
	\begin{array}{ll}
		1+\frac{1}{\delta^{\beta+1}} & \text{for $\beta\neq -1$}\\\\
		1+\log\delta & \text{for $\beta= -1$.}
	\end{array}\right.
\end{align*}
As $m+\beta>0$ by assumption, we conclude that 
\[\lim_{\delta\downarrow0}\delta^m\int_{-T}^T \rho(t,\delta)\, dt =0\]
for all cases under consideration. As noted above, this implies that the 
integrals in \eq{part1} and \eq{part2} tend to zero as $\delta\downarrow0$.
This concludes the proof that $(\rho,u)$ satisfies the weak form
\eq{radial_mass_weak}-\eq{radial_mom_weak} of the radial, isothermal 
Euler system.

We summarize our findings in the following theorem. We recall that the
kink-solution $U_\kaa(\xi)$ refers to the unique solution of the similarity 
ODE \eq{U_ode} on $(-\infty,\xi_w)$ which approaches the critical point $(\xi_\w,U_\w)$ 
with slope $1-\lambda_+$, where $\lambda_+$ is given by \eq{lambdas}.
We also recall the assumption that its limiting value $U^*$ at $\xi=-\infty$ 
is strictly negative (the analysis in Section \ref{constr_conv_div_solns}
shows that this is a non-vacuous assumption). 
%%%%%%%%%%%%%%%%%%%%%%%%%%%%%%%%
\begin{theorem}\label{main_result}
	Consider the radial, isothermal Euler system 
	\eq{mass}-\eq{momentum} with pressure
	function $p=a^2\rho$ in $n=2$ or $3$ space dimensions. With 
	$m=n-1$, choose any $\beta\in(-m,0)$ so that the limiting value
	$U^*$ of the kink-solution $U_\kaa(\xi)$ at $\xi=-\infty$ satisfies
	$U^*<0$. Then, the functions $U(\xi)$ and $\Omega(\xi)$ constructed
	in Section \ref{constr_conv_div_solns} yield, via \eq{fncl_relns},
	a radial weak solution $(\rho(t,r),u(t,r))$ to \eq{mass}-\eq{momentum},  
	according to Definition \ref{rad_symm_weak_soln}. 
	
	In particular, any 
	such solution provides a weak solution $\rho(t,{\bf x}):=\rho(t,|{\bf x}|)$,
	${\bf u}(t,{\bf x}):=u(t,|{\bf x}|)\frac{{\bf x}}{|{\bf x}|}$ to the original, multi-d 
	isothermal system \eq{mass_m_d_isthrml_eul}-\eq{mom_m_d_isthrml_eul}, 
	according to Definition \ref{weak_soln}. Finally, any  
	such solution involves a continuous, focusing wave, followed by an expanding 
	shock wave, and suffers amplitude blowup of its density field at the 
	origin $(t,{\bf x})=(0,0)$, with $\rho(0,{\bf x})\propto |{\bf x}|^\beta$, 
	while its velocity field remains globally bounded.
\end{theorem}
\section{Final remarks}\label{final_rmks}
%%%%%%%%%%%%%%%%%%%%%%%%%%%%%%%%
First, for any fixed time $t$, as $r\to\infty$ the radial speed 
$u(t,r)$ tends to $U^*<0$, while the density $\rho(t,r)$ tends to zero.
However, the latter decay is too slow to give  
bounded total mass. In fact, the solutions constructed above have both
unbounded total mass and unbounded total energy. E.g., the mass density $\rho(t,r)r^m$ grows like 
$r^{\beta +m}$ for $t$ fixed as $r\to\infty$, and the standing assumption that $\beta+m>0$
yields unbounded mass.
A similar calculation shows that the total energy density 
\[E(t,r):=\big[ \textstyle \frac{1}{2}\rho(t,r) u(t,r)^2
+a^2\rho(t,r)\log \rho(t,r)\big]r^m,\]
has unbounded integral at all times. On the other hand, as verified above, mass and energy 
are both locally integrable with respect to space at any fixed times. 

Next, consider the behavior of characteristics $\dot r=u\pm a$
and particle trajectories $\dot r=u$ in the constructed solutions. 
We first note that the only possibility for the path $\xi=\bar\xi$ (constant) 
to be a characteristic, is for $\bar \xi$ to have the value $\xi_\w$. This
yields the ``critical,'' converging 1-characteristic through the origin.
All 1-characteristics below the critical one end up along $\{r=0\}$ 
at negative times (with speed $-a$), while all 1-characteristics 
above it cross $\{t=0\}$ (all with speed $U^*-a$ and at strictly positive distances 
to the origin), and subsequently disappear
into the reflected shock wave propagating along $r=\xi_\s t$. 

Next, all
particle trajectories cross the critical characteristic from below (in the 
$(r,t)$-plane) and proceed to cross $\{t=0\}$ with speed $U^*$.
It follows that there is no ``accumulation'' of particles at the center of 
motion; in particular, the trivial particle trajectory $r(t)\equiv 0$ is
the unique one passing through the origin. Consequently, the density 
$\rho(t,r)$ does not ``contain a Dirac delta'' at time of collapse.
(Solutions of ``cumulative'' type where all, or part, of the mass 
concentrates at the origin at some instance have been considered 
in \cites{kell,am}.)

Finally, let $\{r=\mathfrak c(t)\}$ be any 1-characteristic above the 
critical 1-characteristic $\{\xi=\xi_\w\}$; then $\mathfrak c(0)>0$. 
We could now replace the constructed similarity solution on 
$\{r>\mathfrak c(t)\}$ with a solution (e.g., a simple wave with the same
values along $\{r=\mathfrak c(t)\}$) of 
finite mass and energy in this outer region, without affecting the 
behavior of the solution within $\{r<\mathfrak c(t)\}$. This shows that 
the type of amplitude blowup exhibited by the original similarity solution, 
is possible also in solutions with finite mass and energy.

\bigskip
\paragraph{\bf Acknowledgment:}
This work was supported in part by NSF awards DMS-1813283 (Jenssen) 
and DMS-1714912 (Tsikkou).

%BIBLIOGRAPHY

\end{document}